\documentclass[a4paper,12pt]{article}

\usepackage{amsmath}
\usepackage{amssymb}
\usepackage[OT4]{fontenc}
\usepackage[utf8]{inputenc}
\usepackage[numbers]{natbib}
\usepackage{tikz}
\usepackage{url}
\usepackage[unicode,debug,colorlinks=true,linkcolor=black,citecolor=black,urlcolor=black,pdfusetitle=true,breaklinks=true]{hyperref}

\allowdisplaybreaks

\newcommand{\abs}[1]{\left| #1 \right|}
\newcommand{\okra}[1]{\left( #1 \right)}
\newcommand{\kwad}[1]{\left[ #1 \right]}
\newcommand{\klam}[1]{\left\{ #1 \right\}}
\newcommand{\boole}[1]{{\bf 1}{\klam{#1}}}
\DeclareMathOperator{\id}{id}

\newtheorem{definition}{Definition}
\newtheorem{example}{Example}

\newtheorem{theorem}{Theorem}
\newtheorem{proposition}{Proposition}
\newenvironment*{proof}{\begin{trivlist}\item[]
\noindent\textbf{Proof:}}{$\Box$\par\end{trivlist}}
\newenvironment*{proof*}[1]{\begin{trivlist}\item[]
\noindent\textbf{Proof of #1:}}{$\Box$\par\end{trivlist}}

\author{{\L}ukasz D\k{e}bowski
  \\
  {\normalsize
      Institute of Computer Science, 
      Polish Academy of Sciences
    }
  \\
  {\normalsize
      ul. Jana Kazimierza 5,
      01-248 Warszawa, 
      Poland 
    }
  \\
  {\normalsize e-mail: \url{ldebowsk@ipipan.waw.pl}}
}

\title{From Letters to Words and Back:\\
  Invertible Coding of Stationary Measures}

\date{}

\begin{document}

\begin{titlepage}
\maketitle

\begin{abstract}
  Motivated by problems of statistical language modeling, we consider
  probability measures on infinite sequences over two countable
  alphabets of a~different cardinality, such as letters and words. We
  introduce an invertible mapping between such measures, called the
  normalized transport, that preserves both stationarity and
  ergodicity. The normalized transport applies so called self-avoiding
  codes that generalize comma-separated codes and specialize bijective
  stationary codes. The normalized transport is also connected to the
  usual measure transport via underlying asymptotically mean
  stationary measures. It preserves the ergodic decomposition. The
  normalized transport and self-avoiding codes arise in the problem of
  successive recurrence times. In particular, we show that successive
  recurrence times are ergodic for an ergodic measure, which
  strengthens a~result by Chen Moy from 1959. We also relate the
  entropy rates of processes linked by the normalized transport.
  \\[1em]
  \\[1em]
\end{abstract}

\end{titlepage}

\section{Introduction}
\label{secIntroduction}

Suppose that we take a~text in English and we split it in two
different ways: on each character in one approach, or on spaces only
according to the second method. What is the relationship between
a~stationary model of text where the variables are words and an
analogous stationary model where the variables are letters and spaces?
In the first formal approximation, we may consider a~finite alphabet
of letters $\mathbb{Y}$, a~countably infinite set of words
$\mathbb{X}$, and a~variable-length code
$f:\mathbb{X}\to\mathbb{Y}^*$, where
$\mathbb{Y}^*:=\bigcup_{n\ge 0} \mathbb{Y}^n$ and
$\mathbb{Y}^0:=\klam{\lambda}$ for an empty string $\lambda$.  We want
to ask and answer the question whether there exists a~natural
bijection between stationary ergodic measures over alphabets
$\mathbb{X}$ and $\mathbb{Y}$ that is driven by a~relatively general
code $f$. If so, how can it be expressed formally? Is there a~simple
expression for this bijection that avoids explicit reference to the
frequency interpretation of probabilities?

Answering such a~question may help for applications in information
theory \cite[Example 6]{GrayKieffer80}, for mathematical foundations
of statistical language modeling and relating various statistical laws
of quantitative linguistics \cite[Chapters 10 and 11]{Debowski21}, as
well as for generalization of the ergodic theorem to random
subsequences \cite{ChenMoy59}.  As for related questions, works
\cite{Wolff80, DeMarcken96} studied the computational linguistic
problem of inferring the partition of a~stream of letters (without
delimiting spaces) into a~sequence of contiguous words. This problem
gave rise to a class of grammar-based codes in universal lossless data
compression \cite{KiefferYang00, CharikarOthers05}. A~similar scheme
called byte-pair encoding is presently used for deep learning of large
language models \cite{SenrichHaddowBirch16}.  The relationship between
two (non-stationary) probability distributions over finite sequences
of words and letters was recently studied from a~computational
perspective by \cite{VieiraOthers24}.  The specific problem of
variable-length coding of stationary measures over infinite sequences
of words and letters drew attention of a~few researchers
\cite{CariolaroPierobon77, GrayKieffer80,
  TimoBlackmoreHanlen2007}. See also sources \cite{Cox62, Elga00} for
less related problems in other domains.

In the problem of variable-length coding of stationary measures, some
moderately general methods were developed by us in \cite{Debowski10,
  Debowski21}. The motivating need was to encode Santa Fe
processes. The Santa Fe processes are some simplistic non-ergodic
statistical models of language over a~countably infinite alphabet
\cite{Debowski11b, Debowski21b, Debowski23} that combine Zipf's law
\cite{Zipf35} with non-ergodicity. They provide a baseline theoretical
model of the neural scaling law in deep learning \cite{KaplanOther20,
  Hutter21}, also known as Hilberg's law \cite{Hilberg90,
  CrutchfieldFeldman03}. We wanted to encode Santa Fe processes into
the ternary alphabet while preserving their stationarity and ergodic
decomposition. We also pursued simple expressions for this
transformation since we wanted to explicitly connect the entropy rates
and the power-law rates of mutual information, called Hilberg
exponents \cite{Debowski15d}, of the related stationary measures, in
particular.

Despite apparently broad methods developed in paper \cite{Debowski10},
later resumed and generalized in \cite[Chapters 10]{Debowski21}, the
basic problem of specifying a~simple bijection between stationary
measures over alphabets $\mathbb{X}$ and $\mathbb{Y}$ remained
unstated and unsolved. We did not expect that there can be a~simple
solution within the reach of our hands. It was only coming across an
article by Chen Moy \cite{ChenMoy59}, who demonstrated stationarity of
successive recurrence times for a~stationary ergodic measure, that
made us realize a~clear simple construction after a~period of thinking
in the background.

The present paper reports our findings.  Our construction is based on
so called self-avoiding codes, which are codes over infinite sequences
that satisfy a certain synchronizability condition. On one hand, these
codes generalize extensions of comma-separated codes $f(x)=w(x)c$,
where $w(x)\in(\mathbb{Y}\setminus\klam{c})^*$. On the other hand,
they are isomorphic with a~subclass of bijective stationary codes in
the sense of Gray \cite{Gray90}, to be called chunked codes.  Using
a~self-avoiding code, we introduce a~bijection between measures over
alphabets $\mathbb{X}$ and $\mathbb{Y}$ that preserves both
stationarity and ergodicity.  This one-to-one mapping is called the
normalized transport.  It is connected to the usual measure transport
via underlying asymptotically mean stationary measures. The normalized
transport preserves the ergodic decomposition. The normalized
transport arises, for example, in the problem of successive recurrence
times. Specifically, we will prove that successive recurrence times
are ergodic for an ergodic measure, which extends the result of Chen
Moy \cite{ChenMoy59}.

Another sort of our findings concerns extending the results of
\cite{TimoBlackmoreHanlen2007, Debowski10, Debowski21} that link
certain information measures of stationary processes connected by the
normalized transport. The entropy rate is the basic measure in this
context. Let us recall that any bijective stationary code preserves
the entropy rate of the related stationary measures \cite[Corollary
4.2.5]{Gray90}. Since any self-avoiding code is a composition of a
string-valued bijective stationary code and the infinite
concatenation, we can prove that the entropy rates of stationary
ergodic measures that are linked by the normalized transport are
proportional. The proportionality factor is the asymptotic ratio of
the lengths of corresponding sequences of words and letters.

The contents of this article are as follows. In Section
\ref{secPreliminaries}, we describe our notation. In Section
\ref{secCodes}, we introduce self-avoiding codes, covering both theory
and some examples. In Section \ref{secTransport}, we demonstrate how
the normalized transport operates, stating and proving the main
theorems of this paper. In Section \ref{secRecurrence}, we prove that
successive recurrence times are stationary ergodic for a stationary
ergodic process. In Section \ref{secEntropy}, we discuss how the
entropy rate is transformed by the normalized transport.  Section
\ref{secProblems} sketches some open problems and concludes the paper.

\section{Preliminaries}
\label{secPreliminaries}

To establish the notation, we state the following. Consider measurable
spaces of two-sided infinite sequences
$(\mathbb{X}^{\mathbb{Z}},\mathcal{X}^{\mathbb{Z}})$ and
$(\mathbb{Y}^{\mathbb{Z}},\mathcal{Y}^{\mathbb{Z}})$ over two
countable alphabets $\mathbb{X}$ and $\mathbb{Y}$. For symbols
$x_i\in\mathbb{X}$ and $y_i\in\mathbb{Y}$, we write individual
two-sided sequences as
\begin{align}
  x^{\mathbb{Z}}&:=\ldots x_{-1}x_{0};x_1x_2\ldots \,,
  &
  y^{\mathbb{Z}}&:=\ldots y_{-1}y_{0};y_1y_2\ldots \,.
\end{align}
Mind the semicolon between the $0$-th and the first element. We define
random variables
\begin{align}
  \label{Projections}
  X_i(x^{\mathbb{Z}})&:=x_i,
  &
  Y_i(y^{\mathbb{Z}})&:=y_i.
\end{align}
We have two shift operations
\begin{align}
  T_X x^{\mathbb{Z}}&:=\ldots x_{0}x_{1};x_2x_3\ldots \,,
  &
  T_Y y^{\mathbb{Z}}&:=\ldots y_{0}y_{1};y_2y_3\ldots \,.
\end{align}

Let us choose some sets $\Omega_X\in \mathcal{J}_X$ and
$\Omega_Y\in \mathcal{J}_Y$ that are shift-invariant, i.e., satisfy
$T_X^{-1}\Omega_X=\Omega_X$ and $T_Y^{-1}\Omega_Y=\Omega_Y$. We define
the restricted $\sigma$-fields
\begin{align}
  \mathcal{J}_X&:=\klam{A\cap \Omega_X:A\in\mathcal{X}^{\mathbb{Z}}},
  &
  \mathcal{J}_Y&:=\klam{B\cap \Omega_Y:B\in\mathcal{Y}^{\mathbb{Z}}}.
\end{align}
We also consider two invariant $\sigma$-fields
\begin{align}
  \mathcal{I}_X&:=\klam{A\in\mathcal{J}_X:T_X^{-1}A=A},
  &
  \mathcal{I}_Y&:=\klam{B\in\mathcal{J}_Y:T_Y^{-1}B=B}.
\end{align}
Consider next two measures $P_X$ on $(\Omega_X,\mathcal{J}_X)$ and
$P_Y$ on $(\Omega_Y,\mathcal{J}_Y)$.  Let $Z\in\klam{X,Y}$. A~measure
$P_Z$ is called:
\begin{itemize}
\item asymptotically mean stationary (AMS) if limits
  \begin{align}
    \bar P_Z(C)=
    \lim_{n\to\infty}\frac{1}{n}\sum_{i=0}^{n-1}P_Z(T_Z^{-i} C)
  \end{align}
  exists for all $C\in \mathcal{J}_Z$ (measure $\bar P_Z$
  is called the stationary mean of $P_Z$);
\item stationary if $P_Z(C)=P_Z(T_Z^{-1} C)$
  for all $C\in \mathcal{J}_Z$;  
\item ergodic if $P_Z(C)\in\klam{0,1}$ for all $C\in\mathcal{I}_Z$.
\end{itemize}
Let $\mathbb{A}_Z$, $\mathbb{S}_Z$, and $\mathbb{E}_Z$ be the classes
of AMS, stationary, and stationary ergodic probability measures on
$(\Omega_Z,\mathcal{J}_Z)$, respectively.

Throughout this paper, we adopt $\inf\emptyset:=\infty$ and
$\sup\emptyset:=-\infty$.

\section{Self-avoiding codes}
\label{secCodes}

Let us reconsider the example of a~text in English mentioned in the
first sentence of this article.  The mapping between script or speech
and the stream of contiguous words might be considered a trivial
comma-separated code. In the academic written English, there are
visible spaces between words and also in speech, we have an impression
of hearing intermittent silences. Nonetheless, the last impression is
an illusion \cite{Jelinek97}, whereas the mapping between script and
words is not so simple, as the following example attests.
\begin{example}[English ``not'']
  \label{exNot}
  We present a simple example of context-sensitive separator deletion.
  Consider the mapping between words and letters given by the rules of
  the informal written English. Word \textbf{not} may be spelt as
  \textbf{n't} or \textbf{'t} (after \textbf{can}). Whenever it
  happens, the separator (space) between this contraction and the
  preceding word is deleted. The separator omission does not hinder
  correct parsing of the text into words.
\end{example}

Similar and even more messy examples can be multiplied in other
written languages: definite articles in French, agglutinants in
Polish, non-concatenative morphology in Arabic, lack of spaces between
words in Chinese, word separators implied by the change of script in
Japanese, etc.  Although the mapping between script or speech and the
stream of contiguous words can be sometimes surprisingly complex
across different human languages, their native speakers are usually
able to correctly identify words in the text written or spoken in
their own language.  Despite complexity of such mappings, we expect
that the script and the stream of contiguous words convey the same
information and have similar probabilistic properties.

Inspired by these linguistic examples, we would like to construct an
invertible mapping from stationary ergodic probability measures on
$(\Omega_X,\mathcal{J}_X)$ to stationary ergodic probability measures
on $(\Omega_Y,\mathcal{J}_Y)$.  This will be done through a specific
class of codes, called self-avoiding codes.  On our way, we will
discuss a few other classes of codes, as well.  As the narration
progresses, we will see that self-avoiding codes remain in
a~one-to-one correspondence with a~subclass of bijective stationary
codes in the sense of Gray \cite{Gray90}. This subclass is called
chunked codes. We will also see a few examples of simple self-avoiding
codes motivated by toy statistical language modeling. We suppose that
self-avoiding codes can be more complicated in general.

To provide a road map, the classes of codes to be discussed have been
depicted in terms of Venn diagrams in Figures \ref{figWordToLetters},
\ref{figWordsToLetters}, and \ref{figWordsToWords}. Figure
\ref{figWordToLetters} collects finitary codes
$\mathbb{X} \to \mathbb{Y}^*$ (one word mapped to a string of
letters), Figure \ref{figWordsToLetters} visualizes infinitary codes
$\Omega_X \to \Omega_Y$ (a sequence of words mapped to a sequence of
letters), whereas Figure \ref{figWordsToWords} shows infinitary codes
$\Omega_X \to \Omega_W$, where $\mathbb{W} = \mathbb{Y}^*$ (a sequence
of words mapped to a sequence of contiguous chunks).

\begin{figure}[t]
  \centering
\begin{tikzpicture}
  \def\thirdcircle{(0,-0.4) ellipse[x radius=4, y radius=2.6]}
  \def\fourthcircle{(0,-0.85) ellipse[x radius=3.2, y radius=2]}
  \def\fifthcircle{(0,-1.5) ellipse[x radius=2.4, y radius=1.2]}
  \draw \thirdcircle;
  \draw \fourthcircle;
  \draw \fifthcircle;
  \node at (0,1.6) {Finitary codes};
  \node at (0,0.2) {Comma-embedded codes};
  \node at (0,-1.7) {Comma-separated codes};
\end{tikzpicture}  
\caption{Finitary codes $\mathbb{X} \to
  \mathbb{Y}^*$. \label{figWordToLetters}}
\end{figure}

\begin{figure}[p]
  \centering
\begin{tikzpicture}
  \def\firstcircle{(0,0) ellipse[x radius=4, y radius=3.2]}
  \def\secondcircle{(0,1.4) ellipse[x radius=4, y radius=3.4]}
  \def\thirdcircle{(0,0.6) ellipse[x radius=3.4, y radius=2.3]}
  \def\fourthcircle{(0,0.9) ellipse[x radius=2.4, y radius=0.7]}
  \def\fifthcircle{(0,1.8) ellipse[x radius=3.4, y radius=1.8]}
  \draw \firstcircle;
  \draw \secondcircle;
  \draw \thirdcircle;
  \draw \fourthcircle;
  \draw \fifthcircle;
  \node at (0,4) {Quasi-periodic codes};
  \node at (0,-2.5) {Synchronizable codes};
  \node at (0,-0.4) {Self-avoiding codes};
  \node at (0,2.1) {Code extensions};
  \node at (0,0.9) {Comma-embedded codes};
\end{tikzpicture}  
\caption{Infinitary codes $\Omega_X \to
  \Omega_Y$. \label{figWordsToLetters}}
\end{figure}

\begin{figure}[p]
  \centering
\begin{tikzpicture}
  \def\firstcircle{(0,-0.7) ellipse[x radius=4, y radius=2]}
  \def\secondcircle{(0,0.7) ellipse[x radius=4, y radius=2]}
  \def\thirdcircle{(0,0) ellipse[x radius=3, y radius=1]}
  \draw \firstcircle;
  \draw \secondcircle;
  \draw \thirdcircle;
  \node at (0,-1.8) {Stationary codes};
  \node at (0,1.8) {Bijective codes};
  \node at (0,-0.3) {(Self-avoiding codes)};
  \node at (0,0.3) {Chunked codes};
\end{tikzpicture}  
\caption{Infinitary codes $\Omega_X \to \Omega_W$, where
    $\mathbb{W} = \mathbb{Y}^*$. \label{figWordsToWords}}
\end{figure}

The formal construction will start with a concept of the quasi-period
that tries to locally relate the shift operations on spaces
$(\Omega_X,\mathcal{J}_X)$ and $(\Omega_Y,\mathcal{J}_Y)$.
\begin{definition}[quasi-period, cf. \mbox{\cite[Definition
    10.3]{Debowski21}}]
  For a~measurable function $g:\Omega_X\to\Omega_Y$, the quasi-period
  is the random variable
  \begin{align}
    \label{QuasiPeriod}
    L_g(x^{\mathbb{Z}}):=
    \inf\klam{j\ge 1: T_Y^j g(x^{\mathbb{Z}})=g(T_X x^{\mathbb{Z}})}.
  \end{align}
\end{definition}
\noindent
The above definition assumes $j\ge 1$ in contrast to $j\ge 0$ in
\cite[Definition 10.3]{Debowski21}. This choice leads to a~simpler
theory.

We will be focused on the following classes of infinitary codes:
\begin{definition}[quasi-periodic code, \mbox{\cite[Definition
    10.3]{Debowski21}}]
  A~measurable function $g:\Omega_X\to\Omega_Y$ is called
  a~quasi-periodic code if for all $x^{\mathbb{Z}}\in\Omega_X$, we
  have
  \begin{align}
    \label{QuasiPeriodic}
    L_g(x^{\mathbb{Z}})<\infty.
  \end{align}
\end{definition}
\begin{definition}[self-avoiding code]
  A~measurable function $g:\Omega_X\to\Omega_Y$ is called
  a~self-avoiding code if it is an injective quasi-periodic code and
  for all $x^{\mathbb{Z}}\in\Omega_X$, we have
  \begin{align}
  \label{SelfAvoiding}
    0<j< L_g(x^{\mathbb{Z}})
    \implies
    T_Y^j g(x^{\mathbb{Z}})\not\in g(\Omega_X).
\end{align}
\end{definition}
In information theory, injective codes are also called non-singular
codes. To recall, $g:\Omega_X\to\Omega_Y$ is called injective if
$a^{\mathbb{Z}}= b^{\mathbb{Z}} \iff g(a^{\mathbb{Z}})=
g(b^{\mathbb{Z}})$.
    
Self-avoiding codes will be the essential building blocks of our
theory. Let us explain properties of these codes in a less formal
way. Conditions (\ref{QuasiPeriodic})--(\ref{SelfAvoiding}) were
distilled from the proofs of theorems in Section \ref{secTransport}.
Appealing to linguistic intuitions, conditions
(\ref{QuasiPeriodic})--(\ref{SelfAvoiding}) capture a general
reasonable mapping between a sequence of contiguous words and its
letter-based representation that applies a well defined left-to-right
writing system.  Condition (\ref{QuasiPeriodic}) states that if we
move in the space of words by one position then we move in the space
of letters by a variable positive integer number of
positions. According to condition (\ref{SelfAvoiding}), we cannot
decipher the text in a wrong way if we move in the space of letters by
a wrong number of positions.

The general class of self-avoiding codes seems to be an invention of
this paper, driven by a few prior specific examples.  These examples,
such as the extensions of comma-separated codes, stay close to naive
intuitions about writing systems of human languages.  Let us recall
these instances. Let $\abs{w}:=n$ be the length of string
$w\in\mathbb{Y}^n$.
\begin{example}[code extension, \mbox{\cite[Section 10.2]{Debowski21}}]
  For a~function $f:\mathbb{X}\to\mathbb{Y}^*$ with $\abs{f(x)}\ge 1$
  for all $x\in\mathbb{X}$, called a~finitary code, we define its
  infinitary extension
  $f^{\mathbb{Z}}:\mathbb{X}^{\mathbb{Z}}\to\mathbb{Y}^{\mathbb{Z}}$
  as
  \begin{align}
    \label{Extension}
    f^{\mathbb{Z}}(x^{\mathbb{Z}})
    &:=
      \ldots  f(x_{-1})f(x_0);f(x_1)f(x_2)\ldots \,.
  \end{align}
  Extension $f^{\mathbb{Z}}$ is a~quasi-periodic code with a
  quasi-period $L_{f^{\mathbb{Z}}}(x^{\mathbb{Z}})\le\abs{f(x_1)}$.
\end{example}

The code extension $f^{\mathbb{Z}}$ need not be a self-avoiding code,
for instance if $f$ is a purely unary code
$f:\mathbb{X}\to\klam{\mathbf{0}}^*$ or, more generally, if $f$ is not
prefix-free.  In the probabilistic context, it might make sense to
consider extensions of codes $f:\mathbb{X}\to\mathbb{Y}^*$ such that
$\abs{f(x)}=0$ for some symbols with probability $0<P_X(X_i=x)<1$.
The next definition clearly excludes such a~possibility to guarantee
that the code extension is self-avoiding.
\begin{example}[comma-separated code, \mbox{\cite[Section 10.2]{Debowski21}}]
  An injection $f:\mathbb{X}\to\mathbb{Y}^*$ is called
  a~comma-separated code if $f(x)=w(x)c$, where
  $w(x)\in(\mathbb{Y}\setminus\klam{c})^*$.  Its extension
  $f^{\mathbb{Z}}$ is a~self-avoiding code with the quasi-period
  $L_{f^{\mathbb{Z}}}(x^{\mathbb{Z}})=\abs{f(x_1)}$. Symbol $c$ is
  called the separator of code $f$.
\end{example}

The above construction can be somewhat generalized while preserving
the self-avoiding property.
\begin{example}[comma-embedded code]
  \label{exEmbedded}
  An injection $f:\mathbb{X}\to\mathbb{Y}^*$ is called
  a~comma-embedded code if $f(x)=w(x)cz(x)$, where
  $w(x)\in(\mathbb{Y}\setminus\klam{c})^*$ and
  $z(x)\in(\mathbb{Y}\setminus\klam{c})^n$ with fixed $n$. Its
  extension $f^{\mathbb{Z}}$ is a~self-avoiding code with the
  quasi-period $L_{f^{\mathbb{Z}}}(x^{\mathbb{Z}})=\abs{f(x_1)}$.
\end{example}
\noindent

For instance, the comma-embedded code with $n=1$ was used for encoding
the Santa Fe process as the Oracle process, in order to guarantee that
the latter is unifilar.
\begin{example}[Santa Fe and Oracle processes,
  \mbox{\cite{Debowski11b,Debowski21b}}]
  A Santa Fe process $(X_i)_{i\in\mathbb{Z}}$ is a sequence of pairs
  $X_i=(K_i,Z_{K_i})$, where $(K_i)_{i\in\mathbb{Z}}$ is a certain
  sequence of random natural numbers $K_i:\Omega\to\mathbb{N}$ and
  $(Z_k)_{k\in\mathbb{N}}$ is a certain sequence of random bits
  $Z_k:\Omega\to\klam{\mathbf{0},\mathbf{1}}$. An Oracle process
  $(Y_j)_{j\in\mathbb{Z}}$ is a certain encoding
  $(Y_j)_{j\in\mathbb{Z}}=f^{\mathbb{Z}}((X_i)_{i\in\mathbb{Z}})$ of
  the Santa Fe process, where
  $f:\mathbb{N}\times\klam{\mathbf{0},\mathbf{1}}
  \to\klam{\mathbf{0},\mathbf{1},\mathbf{2}}^*$ is a comma-embedded
  code of form $f(k,z)=w(k)\mathbf{2}z$ with a non-singular (i.e.,
  injective) code $w:\mathbb{N}\to\klam{\mathbf{0},\mathbf{1}}^*$.
\end{example}

It was shown in \cite{Debowski10,Debowski11b,Debowski21,Debowski21b}
that Santa Fe and Oracle processes have the same ergodic decomposition
and exhibit similar information-theoretic properties, even if we
consider the stationary means of processes $(X_i)_{i\in\mathbb{Z}}$
and $(Y_j)_{j\in\mathbb{Z}}$. If $(Z_k)_{k\in\mathbb{N}}$ is a
sequence of fair-coin flips and $(K_i)_{i\in\mathbb{Z}}$ is
distributed according to Zipf's law, then both processes
$(X_i)_{i\in\mathbb{Z}}$ and $(Y_j)_{j\in\mathbb{Z}}$ are strongly
non-ergodic. They also follow a power-law growth of mutual information
between adjacent blocks of text of an increasing length, providing a
simple model of the neural scaling law \cite{KaplanOther20, Hutter21},
known also as Hilberg's law \cite{Hilberg90, CrutchfieldFeldman03}

Self-avoiding codes can be more complicated. To warm up the
imagination, we note that we can shift the origin of the coordinate
system and obtain other examples of self-avoiding codes.
\begin{example}[shifted code]
  Let $g:\Omega_X\to\Omega_Y$ be a self-avoiding code. Then
  $h=T_Y^m\circ g\circ T_X^n:\Omega_X\to\Omega_Y$ with fixed
  $m,n\in\mathbb{Z}$ is also a~self-avoiding code with the
  quasi-period $L_h(x^{\mathbb{Z}})=L_g(T_X^n x^{\mathbb{Z}})$.
\end{example}

Let us try to characterize self-avoiding codes in a more systematic
way. First, we want to relate them to synchronizable codes, previously
discussed in the literature. For this goal, self-avoiding codes can be
characterized by the following construction.
\begin{definition}[accumulations]
  \label{defiAdjusted}
  For a~measurable function $L:\Omega_X\to\mathbb{N}$, we define
  accumulations $S_k:\Omega_X\to\mathbb{Z}$ as random variables
  \begin{align}
    \label{RecurrenceX}
    S_k
    :=
    \begin{cases}
      S_{k+1}-L\circ T_X^{-k},
      &k\in-\mathbb{N},
      \\
      0,
      &k=0,
      \\
      S_{k-1}+L\circ T_X^{k},
      &k\in\mathbb{N}.
    \end{cases}
  \end{align}
\end{definition}
\begin{definition}[adjusted code]
  A~measurable function $g:\Omega_X\to\Omega_Y$ is called a code
  adjusted to function $L:\Omega_X\to\mathbb{N}$ if for
  all $a^{\mathbb{Z}},b^{\mathbb{Z}}\in\Omega_X$:
  \begin{enumerate}
  \item for all $i\in\mathbb{Z}$, we have
    $T_X^i a^{\mathbb{Z}}= b^{\mathbb{Z}} \iff
    T_Y^{S_i(a^{\mathbb{Z}})} g(a^{\mathbb{Z}})= g(b^{\mathbb{Z}})$;
  \item for all $j\not\in \klam{S_k(a^{\mathbb{Z}}): k\in\mathbb{Z}}$, we have
    $T_Y^j g(a^{\mathbb{Z}})\neq g(b^{\mathbb{Z}})$.
  \end{enumerate}
\end{definition}
\begin{proposition}
  \label{theoIterated}
  A~measurable function $g:\Omega_X\to\Omega_Y$ is a code adjusted to
  a~function $L:\Omega_X\to\mathbb{N}$ if and only if $g$
  is a~self-avoiding code and its quasi-period equals $L_g=L$.
\end{proposition}
\begin{proof}
  For a~self-avoiding code $g:\Omega_X\to\Omega_Y$, let us put
  $L=L_g$.  The $i$-th fold application of definition
  (\ref{QuasiPeriod}) yields
  $T_Y^{S_i(a^{\mathbb{Z}})} g(a^{\mathbb{Z}})=g(T_X^i
  a^{\mathbb{Z}})$. Since $g$ is an injection, this implies condition
  (i) of Definition \ref{defiAdjusted}.  The $i$-th fold application
  of definition (\ref{QuasiPeriod}) and condition (\ref{SelfAvoiding})
  yields $T_Y^j g(a^{\mathbb{Z}})\not\in g(\Omega_X)$ for
  $S_i(a^{\mathbb{Z}})< j<S_{i+1}(a^{\mathbb{Z}})$. Hence we derive
  condition (ii) of Definition \ref{defiAdjusted}. Conversely, if
  conditions (i) with $i\in\klam{0,1}$ and (ii) are satisfied then
  $L_g=L<\infty$, condition (\ref{SelfAvoiding}) holds, and $g$ is
  injective.
\end{proof}

Let us recall the weaker condition of a synchronizable code.
\begin{definition}[synchronizable code, \mbox{\cite[Definition
    10.4]{Debowski21}}]
  A~measurable function $g:\Omega_X\to\Omega_Y$ is called
  a~synchronizable code if it is injective and for all
  $a^{\mathbb{Z}},b^{\mathbb{Z}}\in\Omega_X$, we have
  \begin{align}
    \label{Synchronizable}
    T_X^i a^{\mathbb{Z}}=
    b^{\mathbb{Z}}\text{ for an $i\in\mathbb{Z}$}
    \iff
    T_Y^j g(a^{\mathbb{Z}})=
    g(b^{\mathbb{Z}})\text{ for an $j\in\mathbb{Z}$}.
  \end{align}
\end{definition}
Condition (\ref{Synchronizable}) guarantees synchronization, namely,
two shifted sequences of letters are equal if and only if they
represent sequences of words that are shifted by an integer number of
positions. In contrast to self-avoiding codes that capture only
left-to-right scripts, condition (\ref{Synchronizable}) may describe
both left-to-right and right-to-left writing systems.

Conditions (\ref{QuasiPeriodic})--(\ref{SelfAvoiding}) imply condition
(\ref{Synchronizable}) if the code is injective.
\begin{proposition}
  Any self-avoiding code is a synchronizable code.
\end{proposition}
\begin{proof}
  The claim follows by Proposition \ref{theoIterated} and the simple
  fact that any adjusted code is a synchronizable code.
\end{proof}

Now, we would like to prove a more fundamental characterization of
self-avoiding codes. Namely, each self-avoiding code is the
concatenation of a certain uniquely defined chunked code. The
intuitive meaning of this result is that each self-avoiding code maps
a sequence of words into a well defined sequence of contiguous strings
of letters, which we call chunks. These chunks need not correspond to
individual words, however.  Yet, the chunked codes are a subclass of
bijective stationary codes.
 
Thus, in the following, we consider infinite sequences over the
alphabet of strings $\mathbb{W}:=\mathbb{Y}^*$ and we adopt analogous
notations for the derived concepts as in Section
\ref{secPreliminaries}. We depart from the notion of successive
recurrence times considered by Chen Moy \cite{ChenMoy59}. 
\begin{definition}[recurrence times]
  Consider an event $B\in\mathcal{Y}^{\mathbb{Z}}$. Define successive
  recurrence times as partial functions
  \begin{align}
    \label{RecurrenceY}
    R_k
    :=
    \begin{cases}
      \sup\klam{i\le -1: T_Y^{R_{k+1}+i} Y^{\mathbb{Z}}\in B},
      &k\in-\mathbb{N},
      \\
      0,
      &k=0 \text{ and } Y^{\mathbb{Z}}\in B,
      \\
      \inf\klam{i\ge 1: T_Y^{R_{k-1}+i} Y^{\mathbb{Z}}\in B},
      &k\in\mathbb{N}
    \end{cases}
  \end{align}
  on $(\mathbb{Y}^{\mathbb{Z}},\mathcal{Y}^{\mathbb{Z}})$. Define
  event
  $\Omega^B_Y:=\okra{R_k\in\mathbb{Z} \text{ for all }
    k\in\mathbb{Z}}$---on which all these functions are defined.
\end{definition}

These successive recurrence times define a unique parsing of the
sequence of letters into a sequence of chunks.
\begin{definition}[recurrence parsing]
  Let $\mathbb{W}:=\mathbb{Y}^*$ be the alphabet of strings.  The
  recurrence parsing induced by event an
  $B\in\mathcal{Y}^{\mathbb{Z}}$ is a~mapping
  $h:\Omega^B_Y\to\mathbb{W}^{\mathbb{Z}}$ such that
  $h(y^{\mathbb{Z}})=w^{\mathbb{Z}}$ where
  \begin{align}
    w_k:=y_{R_{k-1}(y^{\mathbb{Z}})+1}y_{R_{k-1}(y^{\mathbb{Z}})+2}
    \ldots y_{R_{k}(y^{\mathbb{Z}})}.
  \end{align}
  Strings $w_k$ are called chunks.  We denote the image
  $\Omega_W:=h(\Omega^B_Y)$ and the set
  $\Omega_Y:=\bigcup_{i\in\mathbb{Z}}
  T_Y^i\Omega^B_Y\in\mathcal{I}_Y$.
\end{definition}

There is a natural example of recurrence times for self-avoiding
codes. In fact, we had already given it a different name.
\begin{example}[recurrence of the code image]
  Accumulations $S_k$ of the quasi-period of a self-avoiding code $g$
  are equal to recurrence times $R_k$ of event
  $(Y^{\mathbb{Z}}\in g(\Omega_X))$. Moreover, both are equal to the
  lengths of corresponding chunks $w_k$. Precisely, we have
\begin{align}
  \abs{w_k}=R_k(y^{\mathbb{Z}})=S_k(x^{\mathbb{Z}}),
\end{align}
where we follow definition (\ref{RecurrenceX}) with $L=L_g$ and
definition (\ref{RecurrenceY}) with
$B:=(Y^{\mathbb{Z}}\in g(\Omega_X))$,
$y^{\mathbb{Z}}=g(x^{\mathbb{Z}})$, and
$w^{\mathbb{Z}}=h(y^{\mathbb{Z}})$, where $h$ is the recurrence
parsing induced by event $B$.
\end{example}

Let us denote the identity code
$\id:\mathbb{W}\ni w\mapsto w\in\mathbb{W}$. In our notation, its
extension $\id^{\mathbb{Z}}|_{\Omega_W}:\Omega_W\to\Omega_Y$ is the
mapping that concatenates an infinite sequence of strings into an
infinite sequence of symbols. We observe that each recurrence parsing
is a one-to-one mapping and its inverse mapping is a self-avoiding
code.
\begin{proposition}
  \label{theoRecurrenceInverse}
  Consider an event $B\in\mathcal{Y}^{\mathbb{Z}}$.
  Recurrence parsing $h:\Omega^B_Y\to\Omega_W$ is a~measurable
  bijection with $\Omega_W\in\mathcal{I}_W$ and
  $h^{-1}=\id^{\mathbb{Z}}|_{\Omega_W}:\Omega_W\to\Omega_Y$ is
  a~self-avoiding code.
\end{proposition}
\begin{proof}
By the definition of recurrence parsing $h$, if
$w^{\mathbb{Z}}=h(y^{\mathbb{Z}})$ then
$y^{\mathbb{Z}}=\id^{\mathbb{Z}}(w^{\mathbb{Z}})$.  Suppose that
$h(a^{\mathbb{Z}})=h(b^{\mathbb{Z}})$. Then
\begin{align}
  a^{\mathbb{Z}}
  =
  \id^{\mathbb{Z}}(h(a^{\mathbb{Z}}))=\id^{\mathbb{Z}}(h(b^{\mathbb{Z}}))
  =
  b^{\mathbb{Z}}
  .
\end{align}
Thus function $h:\Omega^B_Y\to\Omega_W=h(\Omega^B_Y)$ is
a~bijection and $h^{-1}=\id^{\mathbb{Z}}|_{\Omega_W}$. Measurability
of function $h$ follows by measurability of preimages of the
appropriate cylinder sets.  We also have $\Omega_W\in\mathcal{I}_W$
since
\begin{align}
  \Omega_W
  =
  h\okra{R_k\in\mathbb{Z} \text{ for all } k\in\mathbb{Z}}
  =
  \bigcap_{i\in\mathbb{Z}} \okra{\id^{\mathbb{Z}}(T_W^i W^{\mathbb{Z}})\in B}
  .
\end{align}

In the following, we show that the inverse
$h^{-1}=\id^{\mathbb{Z}}|_{\Omega_W}$ is a~self-avoiding code.  First
of all, function $h^{-1}:\Omega_W\to\Omega^B_Y\subset\Omega_Y$ is
a~bijection since $h:\Omega^B_Y\to\Omega_W$ is its inverse.  Second,
for $w^{\mathbb{Z}}=h(y^{\mathbb{Z}})$, we observe
\begin{align}
  \label{LR}
  L_{h^{-1}}(w^{\mathbb{Z}})=\abs{w_1}=R_1(y^{\mathbb{Z}}).
\end{align}
In particular, we have $R_1(y^{\mathbb{Z}})<\infty$ for all
$y^{\mathbb{Z}}\in\Omega_Y$ and
$w^{\mathbb{Z}}=h(y^{\mathbb{Z}})$. Thus function $h^{-1}$ satisfies
condition (\ref{QuasiPeriodic}).  Third, suppose that
$a^{\mathbb{Z}},b^{\mathbb{Z}}\in\Omega^B_Y$.  Condition
$T_Y^i a^{\mathbb{Z}}=b^{\mathbb{Z}}$ for an $i\in\mathbb{Z}$ is
equivalent to $T_W^j h(a^{\mathbb{Z}})=h(b^{\mathbb{Z}})$ for an
$j\in\mathbb{Z}$.  Thus function $h^{-1}$ satisfies condition
(\ref{Synchronizable}).  Finally, $T_Y^j y^{\mathbb{Z}}\in \Omega^B_Y$
for $y^{\mathbb{Z}}\in\Omega^B_Y$ implies $j\le 0$ or
$j\ge R_1(y^{\mathbb{Z}})$.  Thus function $h^{-1}$ satisfies
condition (\ref{SelfAvoiding}).
\end{proof}

Proposition \ref{theoRecurrenceInverse} allows to partly characterize
self-avoiding codes using the concept of bijective stationary
codes. The concept of stationary codes is well known in the
literature.
\begin{definition}[stationary code, \mbox{\cite[Section 2.2]{Gray90}}]
  A~measurable function $s:\Omega_X\to\Omega_W$ is called a~stationary
  code if
\begin{align}
  \label{Stationary}
  T_W\circ s=s\circ T_X.
\end{align}
\end{definition}

We also notice this simple fact about bijective stationary codes.
\begin{proposition}
  \label{theoStationaryInverse}
  The inverse of a~bijective stationary code is stationary.
\end{proposition}
\begin{proof}
  For (\ref{Stationary}), we may write
  \begin{align}
    s^{-1}\circ T_W=
    s^{-1}\circ T_W\circ s\circ s^{-1}=s^{-1}\circ s\circ T_X\circ s^{-1}=
    T_X\circ s^{-1}.
  \end{align}
  Hence the stationarity of $s^{-1}:\Omega_W\to\Omega_X$ follows.
\end{proof}

Stationary codes may operate both symbolwise and non-locally, in a
sense.  As noticed in by Gray \cite[Section 2.2]{Gray90}, any
stationary code $s$ can be written as
$s(x^{\mathbb{Z}})=w^{\mathbb{Z}}$ where
\begin{align}
  w_k=s_0(T_X^k x^{\mathbb{Z}})
\end{align}
for a~certain function $s_0:\Omega_X\to\mathbb{W}$.  By Proposition
\ref{theoStationaryInverse}, for a bijective stationary code $s$, we
may write its inverse as $s^{-1}(w^{\mathbb{Z}})=x^{\mathbb{Z}}$ where
\begin{align}
  x_k=s^{-1}_0(T_W^k w^{\mathbb{Z}})
\end{align}
for a~certain function $s^{-1}_0:\Omega_W\to\mathbb{X}$.

In fact, each self-avoiding code results from concatenating a~certain
bijective stationary code. This corresponding bijective stationary
code is called the chunked code.
\begin{definition}[chunked code]
  Let function $g:\Omega_X\to\Omega_Y$ be a~self-avoiding code. We
  define event $B:=(Y^{\mathbb{Z}}\in g(\Omega_X))$ and consider
  function $h:\Omega_Y\to\Omega_W$ being the recurrence parsing
  induced by event $B$. Composition
  \begin{align}
    s:=h\circ g
  \end{align}
  is called the chunked code induced by code $g$.
\end{definition}

Now comes the aforementioned main result of this section. 
\begin{proposition}
  \label{theoChunkedCode}
  Any chunked code $s:\Omega_X\to\Omega_W$ is a~bijective stationary
  code.  Moreover, having function $s$, we may recover the original
  self-avoiding code $g$ through composition
\begin{align}
  \label{SelfAvoidingStationary}
  g=\id^{\mathbb{Z}}\circ s.
\end{align}
\end{proposition}
\begin{proof}
  Condition (\ref{Stationary}) follows by observation
  $\abs{w_1}=R_1(y^{\mathbb{Z}})=L_g(x^{\mathbb{Z}})$ for
  $B:=(Y^{\mathbb{Z}}\in g(\Omega_X))$,
  $y^{\mathbb{Z}}=g(x^{\mathbb{Z}})$, and
  $w^{\mathbb{Z}}=h(y^{\mathbb{Z}})$.  Equality
  (\ref{SelfAvoidingStationary}) follows by
  $h^{-1}=\id^{\mathbb{Z}}|_{\Omega_W}$.
\end{proof}

Thus the class of self-avoiding codes is isomorphic with the class of
chunked codes---a~subclass of bijective stationary codes. However, not
every bijective stationary code yields a~self-avoiding code through
concatenation.
\begin{example}[purely unary code]
  For $\mathbb{X}=\mathbb{N}$, $\mathbb{W}=\klam{\mathbf{0}}^*$, and
  $s_0(x^{\mathbb{Z}})=\mathbf{0}^{x_0}$, the respective stationary
  code $s$ is a~bijection with $s^{-1}_0(w^{\mathbb{Z}})=\abs{w_0}$
  but function $g=\id^{\mathbb{Z}}\circ s$ is not an injection since
  it maps every input sequence to constant sequence
  $\mathbf{0}^{\mathbb{Z}}$.
\end{example}

Another twist, the linear ordering of output chunks need not be the
linear ordering of input words.
\begin{example}[swapped code]
  Let $\mathbb{X}=\klam{\mathbf{1},\mathbf{2},\mathbf{3}}$ and
  $\mathbb{Y}=\klam{\mathbf{0},\mathbf{1},\mathbf{2},\mathbf{3}}$.  We
  define code
  $f:\klam{\mathbf{1}}^*\times\klam{\mathbf{2},\mathbf{3}}\to
  \klam{\mathbf{0},\mathbf{1},\mathbf{2},\mathbf{3}}^*$ as
  $f(\mathbf{1}^n\mathbf{2}):=\mathbf{2}\mathbf{0}^n$ and
  $f(\mathbf{1}^n\mathbf{3}):=\mathbf{1}^n\mathbf{3}$. The extension
  of code $f$ is a self-avoiding code with recurrence times
  $R_k=k$. Hence the chunked code of this self-avoiding code consists
  of chunks of length $1$ but is not an extension of a code
  $\mathbb{X}\to\mathbb{Y}$.
\end{example}

Hence the class of self-avoiding codes seems large and hard to
enumerate by examples. As we have mentioned, we suppose that
self-avoiding codes entail all reasonable writing systems of natural
languages.

\section{Normalized transport}
\label{secTransport}

Having defined various types of infinitary codes, let us inspect the
problem of measure transport.  Let $P_Y$ be a~measure on
$(\Omega_Y,\mathcal{J}_Y)$ and let $P_X$ be a~measure on
$(\Omega_X,\mathcal{J}_X)$.  Let function $g:\Omega_X\to\Omega_Y$ be
measurable.  The usually considered measure transport
\begin{align}
  \label{PlainInverse}
  P_Y(B)=P_X(g^{-1}(B)),
  \quad B\in \mathcal{J}_Y,
\end{align}
yields a~probability measure $P_Y$ for any probability measure $P_X$.
As shown in \cite{Debowski10, Debowski21}, mapping
(\ref{PlainInverse}) does not preserve stationarity (it is not true
that $P_X\in\mathbb{S}_X\implies P_Y\in\mathbb{S}_Y$) but, for
a~quasi-periodic code $g$, it preserves the AMS property
($P_X\in\mathbb{A}_X\implies P_Y\in\mathbb{A}_Y$).  Suppose now that
function $g:\Omega_X\to\Omega_Y$ is injective but not necessarily
surjective. Then for measure transport
\begin{align}
  \label{Plain}
  P_X(A)=P_Y(g(A)),
  \quad A\in \mathcal{J}_X,
\end{align}
measure $P_X$ fails to be a~probability measure for a~probability
measure $P_Y$ when $P_Y(g(\Omega_X))<P_Y(\Omega_Y)$. For example, for
a~comma-separated code $f$ with separator $c$, we have
\begin{align}
  \label{CommaSeparatedNormalization}
  P_Y(f^{\mathbb{Z}}(\Omega_X))\le P_Y(Y_0=c).
\end{align}
However, if we start with a~probability measure $P_X$ then mapping
(\ref{Plain}) is the inverse of mapping (\ref{PlainInverse}). In the
following, mapping (\ref{PlainInverse})--(\ref{Plain}) is called the
plain transport.

As we have said, the plain transport, under some conditions on the
infinitary code $g$, preserves the AMS property but it does not
preserve the stationarity. The stationary means of the corresponding
AMS measures are linked by a mapping that differs to
(\ref{PlainInverse})--(\ref{Plain}). In fact, this
stationarity-preserving mapping is not very complicated but it seems
to work mostly for self-avoiding codes. Inspired by the result of
\cite{ChenMoy59}, we notice that it suffices to normalize construction
(\ref{Plain}) to obtain a~stationary ergodic probability measure on
$(\Omega_X,\mathcal{J}_X)$ from a~stationary ergodic probability
measure on $(\Omega_Y,\mathcal{J}_Y)$. This normalized mapping can be
inverted using previous results of \cite{Debowski10, Debowski21}.

Before we state the appropriate theorem, we need a~few definitions.
\begin{definition}[spread]
  For a~measurable function $g:\Omega_X\to\Omega_Y$, the spreads
  are defined as random variables
  \begin{align}
    \label{Spread}
    G_g(B)(x^{\mathbb{Z}})
    :=\sum_{l=0}^{L_g(x^{\mathbb{Z}})-1}\boole{T_Y^l g(x^{\mathbb{Z}})\in B},
    \quad B\in \mathcal{J}_Y,
  \end{align}
  where random variable $L_g(x^{\mathbb{Z}})$ is the quasi-period.  
\end{definition}
Notice that $0\le G_g(B)\le L_g$, whereas $L_g\ge 1$.

We also distinguish two subclasses of AMS measures
\begin{align}
  \mathbb{A}_X^g&:=\klam{P_X\in\mathbb{A}_X: \int L_g d\bar P_X(\cdot|\mathcal{I}_X)<\infty \text{
  $P_X$-a.s.}},
  \\
  \mathbb{A}_Y^g&:=\klam{P_Y\in\mathbb{A}_Y: \bar P_Y(g(\Omega_X)|\mathcal{I}_Y)>0
  \text{ $P_Y$-a.s.}}.
\end{align}
Let $Z\in\klam{X,Y}$.  Subsequently, we
use
\begin{align}
  \mathbb{S}_Z^g&:=\mathbb{S}_Z\cap \mathbb{A}^g_Z,
  &
  \mathbb{E}_Z^g&:=\mathbb{E}_Z\cap \mathbb{A}^g_Z.
\end{align}
\begin{theorem}
  \label{theoEnsemble} 
  Let $g:\Omega_X\to\Omega_Y$ be a~self-avoiding code.  Consider
  conditions:
  \begin{align}
    \label{EnsembleInverse}
    P_Y(B)
    &=
      \frac{\int G_g(B) dP_X}{\int L_g dP_X},
    \quad
    B\in \mathcal{J}_Y,
    \\
    \label{Ensemble}
    P_X(A)
    &=
      \frac{P_Y(g(A))}{P_Y(g(\Omega_X))},
    \quad
    A\in \mathcal{J}_X.
  \end{align}
  We claim that:
\begin{enumerate}
\item $P_X\in\mathbb{E}_X^g$ and (\ref{EnsembleInverse}) holds if and
  only if $P_Y\in\mathbb{E}_Y^g$ and (\ref{Ensemble}) holds.
\item If conditions mentioned in (i) hold then
  \begin{align}
    \label{Normalization}
    \int L_g dP_X
    =\frac{1}{P_Y(g(\Omega_X))}.
  \end{align}
\end{enumerate}
\end{theorem}
We will prove Theorem \ref{theoEnsemble} in Section \ref{secEnsemble}.

We call bijection (\ref{EnsembleInverse})--(\ref{Ensemble}) the
normalized transport. It can be easily extended to certain non-ergodic
measures.  For this aim, we define a~bijection
$C_g:\mathbb{E}_X^g\to\mathbb{E}_Y^g$ such that $C_g(P_X):=P_Y$ given
by (\ref{EnsembleInverse}) and $C_g^{-1}(P_Y):=P_X$ given by
(\ref{Ensemble}). Mapping $C_g$ is measurable respect to the standard
$\sigma$-fields of stationary ergodic measures $\mathcal{E}_X$ and
$\mathcal{E}_Y$. We denote their croppings as
\begin{align}
  \mathcal{E}_X^g&:=\klam{A\cap\mathbb{A}^g_X:A\in\mathcal{E}_X},
  &
  \mathcal{E}_Y^g&:=\klam{B\cap\mathbb{A}^g_Y:B\in\mathcal{E}_Y}.
\end{align}
By the ergodic decomposition \cite{Gray09},
\cite[Section 4.1]{Debowski21}, for any measures
$P_X\in\mathbb{S}_X^g$ and $P_Y\in\mathbb{S}_Y^g$ there exist unique
prior probability measures $\nu_X$ on
$(\mathbb{E}_X^g,\mathcal{E}_X^g)$ and $\nu_Y$ on
$(\mathbb{E}_Y^g,\mathcal{E}_Y^g)$ such that
\begin{align}
  P_X(A)&:=\int f_X(A) d\nu_X(f_X), 
  &
    P_Y(B)&:=\int f_Y(B) d\nu_Y(f_Y),
\end{align}
for $A\in\mathcal{J}_X$ and $B\in\mathcal{J}_Y$. Moreover, we have
\begin{align}
  \nu_X(A)&=P_X(P_X(\cdot|\mathcal{I}_X)\in A),
  &
  \nu_Y(B)&=P_Y(P_Y(\cdot|\mathcal{I}_Y)\in B),
\end{align}
for  $A\in\mathcal{E}_X^g$ and $B\in\mathcal{E}_Y^g$.
If the prior measures are linked by measure transport
$\nu_Y=\nu_X\circ C_g^{-1}$ or equivalently $\nu_X=\nu_Y\circ C_g$
then we obtain a~measurable bijection of stationary measures
\begin{align}
  P_Y(B)
  &=\int f_Y(B) d\nu_Y(f_Y)=\int C_g (f_X)(B) d\nu_X(f_X)
  \nonumber\\
  &=\int\frac{\int G_g(B) df_X}{\int L_g df_X}d\nu_X(f_X)
    =\int\frac{\int G_g(B) dP_X(\cdot|\mathcal{I}_X)}%
    {\int L_g dP_X(\cdot|\mathcal{I}_X)} dP_X,
          \label{TrajectoryInverse}
  \\
  P_X(A)
  &=\int f_X(A) d\nu_X(f_X)=\int C_g^{-1} (f_Y)(A) d\nu_Y(f_Y)
  \nonumber\\
  &=\int\frac{f_Y(g(A))}{f_Y(g(\Omega_X))}d\nu_Y(f_Y)
    =\int\frac{P_Y(g(A)|\mathcal{I}_Y)}{P_Y(g(\Omega_X)|\mathcal{I}_Y)} dP_Y.
          \label{Trajectory}
\end{align}
In view of the above remarks, it is obvious that
$P_X\in\mathbb{S}_X^g$ and (\ref{TrajectoryInverse}) holds if and only
if $P_Y\in\mathbb{S}_Y^g$ and (\ref{Trajectory}) holds. The
generalized bijection (\ref{TrajectoryInverse})--(\ref{Trajectory}) is
still called the normalized transport.

Since $\nu_Y=\nu_X\circ C_g^{-1}$, the normalized transport
(\ref{TrajectoryInverse})--(\ref{Trajectory}) preserves the ergodic
decomposition of transformed processes.  Moreover, it can be connected
to the plain transport (\ref{PlainInverse})--(\ref{Plain}) for the
respective AMS measures. This connection has to do with the frequency
interpretation of stationary means via the Birkhoff ergodic theorem.
\begin{theorem}
  \label{theoTrajectoryAMS}
  Let $g:\Omega_X\to\Omega_Y$ be a~self-avoiding code.  We claim that:
  \begin{enumerate}
  \item $P_X\in\mathbb{A}_X^g$ and (\ref{PlainInverse}) holds if and
    only if $P_Y\in\mathbb{A}_Y^g$ and (\ref{Plain}) holds with
    equality $P_Y(g(\Omega_X))=1$.
  \item If conditions mentioned in (i) hold then conditions
    (\ref{TrajectoryInverse})--(\ref{Trajectory}) are satisfied with
    substitutions $P_X\to\bar P_X$ and $P_Y\to\bar P_Y$.
  \end{enumerate}
\end{theorem}
We will prove Theorem \ref{theoTrajectoryAMS} in Section
\ref{secTrajectoryAMS}.

As for the precedents of Theorems \ref{theoEnsemble} and
\ref{theoTrajectoryAMS} in the literature, there are five most
relevant propositions, cf. \cite[Example 6]{GrayKieffer80},
\cite{Debowski10, Debowski21, ChenMoy59}. We rewrite them in the
notation of this paper. The first two findings concern quasi-periodic
codes that are not necessarily synchronizable.
\begin{proposition}[\mbox{\cite[Theorem 10.7]{Debowski21}}]
  \label{theoEnsembleWeak} 
  Let $g:\Omega_X\to\Omega_Y$ be a~quasi-periodic code.  If
  $P_X\in\mathbb{S}_X^g$ and (\ref{EnsembleInverse}) holds then
  $P_Y\in\mathbb{S}_Y$.
\end{proposition}
\begin{proposition}[\mbox{\cite[Theorem 10.8]{Debowski21}}]
  \label{theoTrajectoryAMSWeak}
  Let $g:\Omega_X\to\Omega_Y$ be a~quasi-periodic code. If
  $P_X\in\mathbb{A}_X^g$ and (\ref{PlainInverse}) holds then
  $P_Y\in\mathbb{A}_Y$ and stationary mean $\bar P_Y$ is given by
  (\ref{TrajectoryInverse}) with substitutions $P_X\to\bar P_X$ and
  $P_Y\to\bar P_Y$.
\end{proposition}
The next two statements concern synchronizable codes.
\begin{proposition}[\mbox{\cite[Theorem 10.10]{Debowski21}}]
  \label{theoSynchro}
  For a~synchronizable code
  $g:\Omega_X\rightarrow\Omega_Y$, we
  have
  \begin{align}
    \label{SynchroInverse}
    \klam{B\cap g(\Omega_X): B\in \mathcal{I}_Y}
    &=g(\mathcal{I}_X).
  \end{align}
\end{proposition}
\begin{proposition}[\mbox{\cite[Theorem 10.11]{Debowski21}}]
  \label{theoSynchroSpectrum}
  Let $g:\Omega_X\to\Omega_Y$ be a~synchronizable code. Let
  $P_X\in\mathbb{A}_X$ and $P_Y\in\mathbb{A}_Y$.  If
  (\ref{PlainInverse}) holds then $\Sigma(P_X)=\Sigma(P_Y)$ for the
  ergodic spectra
  \begin{align}
   \Sigma(P_X)&:=\klam{P_X(A): A\in\mathcal{I}_X},
   &
   \Sigma(P_Y)&:=\klam{P_Y(B): B\in\mathcal{I}_Y}.
  \end{align}
\end{proposition}

As we can see, the main innovation of the present paper is to find out
the explicit inverse of mappings (\ref{EnsembleInverse}) and
(\ref{TrajectoryInverse}), namely, mappings (\ref{Ensemble}) and
(\ref{Trajectory}). This construction works specifically for self-avoiding
codes, as we have already stressed.

\subsection{Proof of Theorem \ref{theoEnsemble}}
\label{secEnsemble}

Throughout this subsection, we assume a~self-avoiding code $g$.

\subsubsection{Mutual invertibility}

First we show that operation (\ref{Ensemble}) is the inverse of
(\ref{EnsembleInverse}). The conditions for this invertibility are
asymmetric. Assume (\ref{EnsembleInverse}) and an arbitrary
probability measure $P_X$ such that $\int L_g dP_X<\infty$.  By
condition (\ref{SelfAvoiding}) and by $g$ being injective, we obtain
\begin{align}
  \frac{P_Y(g(A))}{P_Y(g(\Omega_X))}
  &=
    \frac{\int G_g(g(A)) dP_X}{\int G_g(g(\Omega_X)) dP_X}
    \nonumber\\
  &=
  \frac{\int
  \sum_{l=0}^{L_g(x^{\mathbb{Z}})-1}\boole{T_Y^l g(x^{\mathbb{Z}})\in g(A)}
  dP_X(x^{\mathbb{Z}})}{\int
  \sum_{l=0}^{L_g(x^{\mathbb{Z}})-1}\boole{T_Y^l g(x^{\mathbb{Z}})\in g(\Omega_X)}
    dP_X(x^{\mathbb{Z}})}
    \nonumber\\
  &=
    \frac{\int
    \boole{g(x^{\mathbb{Z}})\in g(A)}
    dP_X(x^{\mathbb{Z}})}{\int
    \boole{g(x^{\mathbb{Z}})\in g(\Omega_X)}
    dP_X(x^{\mathbb{Z}})}
    =
    \frac{\int
    \boole{x^{\mathbb{Z}}\in A}
    dP_X(x^{\mathbb{Z}})}{\int
    dP_X(x^{\mathbb{Z}})}
    \nonumber\\
  &=
    \frac{P_X(A)}{P_X(\Omega_X)}=P_X(A).
\end{align}
Hence $P_Y(g(\Omega_X))>0$ and operation (\ref{Ensemble}) is the
inverse of (\ref{EnsembleInverse}).

Next, we show that (\ref{EnsembleInverse}) is the inverse of
(\ref{Ensemble}). Assume (\ref{Ensemble}) and a~stationary ergodic
probability measure $P_Y$ such that $P_Y(g(\Omega_X))>0$. By condition
(\ref{SelfAvoiding}), we may write the quasi-period as
\begin{align}
  L_g(x^{\mathbb{Z}})=\min\klam{j\ge 1: T_Y^j g(x^{\mathbb{Z}})\in g(\Omega_X)}
\end{align}
and the spreads as
\begin{align}
  G_g(B)(x^{\mathbb{Z}})
  &
    =\sum_{l=0}^{L_g(x^{\mathbb{Z}})-1}
    \boole{T_Y^l g(x^{\mathbb{Z}})\in B}
    \nonumber\\
  &
    =\sum_{l=0}^{\infty}\boole{T_Y^l g(x^{\mathbb{Z}})\in B}
    \prod_{j=1}^l\boole{T_Y^j g(x^{\mathbb{Z}})\not\in g(\Omega_X)}.
\end{align}
Hence for a~stationary ergodic $P_Y$, we obtain
\begin{align}
  &\int_{g(\Omega_X)} G_g(B)\circ g^{-1} dP_Y
    \nonumber\\
  &=\int
    \kwad{
    \sum_{l=0}^{\infty}\boole{T_Y^l y^{\mathbb{Z}}\in B}
    \boole{y^{\mathbb{Z}}\in g(\Omega_X)}
    \prod_{j=1}^l\boole{T_Y^j y^{\mathbb{Z}}\not\in g(\Omega_X)}
    }
    dP_Y(y^{\mathbb{Z}})
    \nonumber\\
  &=\int
    \boole{y^{\mathbb{Z}}\in B}
    \kwad{
    \sum_{l=0}^{\infty}\boole{T_Y^{-l} y^{\mathbb{Z}}\in g(\Omega_X)}
    \prod_{j=0}^{l-1}\boole{T_Y^{-j} y^{\mathbb{Z}}\not\in g(\Omega_X)}
    }
    dP_Y(y^{\mathbb{Z}})
    \nonumber\\
  &=\int
    \boole{y^{\mathbb{Z}}\in B}
    dP_Y(y^{\mathbb{Z}})
    =P_Y(B)
  \label{InfiniteSpread}
\end{align}
since
\begin{align}
  \int
  \kwad{
  \sum_{l=0}^{\infty}\boole{T_Y^{-l} y^{\mathbb{Z}}\in g(\Omega_X)}
  \prod_{j=0}^{l-1}\boole{T_Y^{-j} y^{\mathbb{Z}}\not\in g(\Omega_X)}
  }
  dP_Y(y^{\mathbb{Z}})
  =
  1
\end{align}
by the Poincar\'e recurrence theorem for $P_Y(g(\Omega_X))>0$.  As
a~result, by (\ref{InfiniteSpread}) and $g$ being injective, we obtain
\begin{align}
  \frac{\int G_g(B) dP_X}{\int L_g dP_X}
  &=
    \frac{\int G_g(B) dP_X}{\int G_g(\Omega_Y) dP_X}
    =
    \frac{\int G_g(B) dP_Y(g(\cdot)|g(\Omega_X))}%
    {\int G_g(\Omega_Y) dP_Y(g(\cdot)|g(\Omega_X))}
   \nonumber\\
  &=
    \frac{\int_{\Omega_X} G_g(B) d(P_Y\circ g)}%
    {\int_{\Omega_X} G_g(\Omega_Y) d(P_Y\circ g)}
    =
    \frac{\int_{g(\Omega_X)} G_g(B)\circ g^{-1} dP_Y}%
    {\int_{g(\Omega_X)} G_g(\Omega_Y)\circ g^{-1} dP_Y}
   \nonumber\\
  &=
    \frac{P_Y(B)}{P_Y(\Omega_Y)}=P_Y(B).
\end{align}
Thus $\int L_g dP_X<\infty$ and operation (\ref{EnsembleInverse}) is
the inverse of (\ref{Ensemble}).

\subsubsection{Conservation of stationarity}

Now we will show that stationarity is preserved.  Suppose that $P_X$
is stationary and assume (\ref{EnsembleInverse}). Stationarity of
$P_Y$ follows by Proposition \ref{theoEnsembleWeak} shown as
\cite[Theorem 10.7]{Debowski21}. We quote the proof for
completeness and comparison. Namely, we observe
\begin{align}
  P_Y(T_Y^{-1} B)-P_Y(B)
  &=
    \frac{\int \kwad{G_g(T_Y^{-1} B) - G_g(B)} dP_X}{\int L_g dP_X}
    \nonumber\\
  &=
  \frac{\int
  \kwad{\boole{T_Y^{L_g} g(\cdot)\in B}-\boole{g(\cdot)\in B}}
    dP_X}{\int L_g dP_X}
    \nonumber\\
  &=
  \frac{\int
  \kwad{\boole{g(T_X(\cdot))\in B}-\boole{g(\cdot)\in B}}
  dP_X}{\int L_g dP_X}
  =
  0.
\end{align}
Hence $P_Y$ is stationary.

Subsequently, suppose that $P_Y$ is stationary and assume
(\ref{Ensemble}).  Sets $A_l:=A\cap(L_g=l)$ for $l\in\klam{0,1,2,\ldots }$
form a~partition of an event $A$.  Similarly, sets $T_X A_l$ for
$l\in\klam{0,1,2,\ldots }$ form a~partition of event $T_X A$.  We have
$g(T_X A_l)=T_Y^l g(A_l)$.  Hence, since $g$ is an injection, we
obtain
\begin{align}
  P_X(T_X A)-P_X(A)
  &=
    \frac{P_Y(g(T_X A))-P_Y(g(A))}{P_Y(g(\Omega_X))}
    \nonumber\\
  &=
    \frac{P_Y(g(T_X \bigcup_{l=0}^\infty A_l))
    -P_Y(g(\bigcup_{l=0}^\infty A_l))}{P_Y(g(\Omega_X))}
    \nonumber\\
  &=
    \frac{\sum_{l=0}^\infty P_Y(g(T_X A_l))
    -\sum_{l=0}^\infty P_Y(g(A_l))}{P_Y(g(\Omega_X))}
    \nonumber\\
  &=
    \frac{\sum_{l=0}^\infty P_Y(T_Y^l g(A_l))
    -\sum_{l=0}^\infty P_Y(g(A_l))}{P_Y(g(\Omega_X))}
    =0.
\end{align}
Hence $P_X$ is stationary.

\subsubsection{Conservation of ergodicity}

Finally, we demonstrate that ergodicity is conserved. Suppose that
$P_X$ is ergodic and assume (\ref{EnsembleInverse}). Let $B$ be
a~$T_Y$-invariant set. By equality (\ref{SynchroInverse}), there
exists a~$T_X$-invariant set $A$ such that $B\cap g(\Omega_X)=g(A)$
and $P_X(A)\in\klam{0,1}$.  Thus
\begin{align}
  P_Y(B)
  &=
    \frac{\int G_g(B) dP_X}{\int L_g dP_X}
    =
    \frac{\int \sum_{l=0}^{L_g-1}\boole{T_Y^l g(\cdot)\in B}
    dP_X}{\int L_g dP_X}
    \nonumber\\
  &=
    \frac{\int L_g\boole{g(\cdot)\in B}
    dP_X}{\int L_g dP_X}
    =
    \frac{\int L_g\boole{g(\cdot)\in B\cap g(\Omega_X)}
    dP_X}{\int L_g dP_X}
    \nonumber\\
  &=
    \frac{\int L_g\boole{g(\cdot)\in g(A)}
    dP_X}{\int L_g dP_X}
    =
    \frac{\int_A L_g
    dP_X}{\int L_g dP_X}=P_X(A).
\end{align}
Hence $P_Y$ is ergodic.

Conversely, suppose that $P_Y$ is ergodic and assume
(\ref{Ensemble}). Let $A$ be a~$T_X$-invariant set. By equality
(\ref{SynchroInverse}), there is a~$T_Y$-invariant set $B$ such that
$B\cap g(\Omega_X)=g(A)$ and $P_Y(B)\in\klam{0,1}$.  Thus
\begin{align}
  P_X(A)
  &=
    \frac{P_Y(g(A))}{P_Y(g(\Omega_X))}
    =
    \frac{P_Y(B\cap g(\Omega_X))}{P_Y(g(\Omega_X))}
    =
    P_Y(B).
\end{align}
Hence $P_X$ is ergodic.

\subsubsection{Equality (\ref{Normalization})}
\label{secNormalization}

Suppose that probability measures $P_X$ and $P_Y$ linked by
(\ref{EnsembleInverse})--(\ref{Ensemble}) are stationary ergodic and
conditions $\int L_g dP_X<\infty$ and $P_Y(g(\Omega_X))>0$ hold. Then
in view of condition (\ref{SelfAvoiding}), we obtain
\begin{align}
  P_Y(g(\Omega_X))=\frac{\int G_g(g(\Omega_X)) dP_X}{\int L_g dP_X}
  =\frac{\int 1 dP_X}{\int L_g dP_X}=\frac{1}{\int L_g dP_X}.
\end{align}

\subsection{Proof of Theorem \ref{theoTrajectoryAMS}}
\label{secTrajectoryAMS}

Throughout this subsection, we assume a~self-avoiding code $g$.

\subsubsection{Mutual invertibility}

Since $g$ is an injection then $P_Y$ is given by (\ref{PlainInverse})
for a~probability measure $P_X$ if and only if $P_X$ is given by
(\ref{Plain}) for a~probability measure $P_Y$ and
$P_Y(g(\Omega_X))=1$.

\subsubsection{Conservation of the AMS property}

We will show that the AMS property is conserved. This will be done by
evaluating the stationary means with the help of the ergodic theorem.
Let $Z\in\klam{X,Y}$. For $C\in\mathcal{J}_Z$, we define the relative
frequencies
\begin{align}
  \Phi_Z(C)(z^{\mathbb{Z}})&:=
  \lim_{n\to\infty}\frac{1}{n}\sum_{i=0}^{n-1} \boole{T_Z^i z^{\mathbb{Z}}\in C}.
\end{align}
As discussed in \cite{Rechard56, GrayKieffer80}, \cite[Theorem
10.4]{Debowski21}, measure $P_Z$ is AMS if and only if limits
$\Phi_Z(C)$ exist $P_Z$-almost surely for all $C\in\mathcal{J}_Z$. As
discussed in \cite{Rechard56, GrayKieffer80}, \cite[Theorem
10.2]{Debowski21}, for an AMS measure $P_Z$, we have
\begin{align}
  \label{StationaryMeanFrequency}
  \bar P_Z(C)=\int \Phi_Z(C) d\bar P_Z=\int \Phi_Z(C) dP_Z.
\end{align}
since measures $P_Z$ and $\bar P_Z$ are equal on the invariant
$\sigma$-field $\mathcal{I}_Z$ and $\Phi_Z(C)$ is
$\mathcal{I}_Z$-measurable.

Suppose now that $P_Y$ is given by (\ref{PlainInverse}) for an AMS
measure $P_X$ with $\int L_g d\bar P_X(\cdot|\mathcal{I}_X)<\infty$
$P_X$-almost surely. Then $\Phi_Y(B)$ exists $P_Y$-almost surely if
$\Phi_Y(B)\circ g$ exists $P_X$-almost surely. The $P_X$-almost sure
existence of $\Phi_Y(B)\circ g$ was established in Proposition
\ref{theoTrajectoryAMSWeak} shown as \cite[Theorem
10.7]{Debowski21}. We quote the proof for completeness and further
reference. Namely, by the Birkhoff ergodic theorem \cite{Birkhoff32,
  Garsia65} and by $g$ being quasi-periodic, we have
\begin{align}
  \Phi_Y(B)(g(x^{\mathbb{Z}}))
  &=
    \lim_{n\to\infty}\frac{1}{n}\sum_{i=0}^{n-1}
    \boole{T_Y^i g(x^{\mathbb{Z}})\in B}
    \nonumber\\
  &=
  \lim_{m\to\infty}\frac{\sum_{j=0}^{m-1} G_g(B)(T_X^j x^{\mathbb{Z}})
    }{\sum_{j=0}^{m-1} L_g(T_X^j x^{\mathbb{Z}})}
    \nonumber\\
  &=
    \frac{\int G_g(B) d\bar P_X(\cdot|\mathcal{I}_X)(x^{\mathbb{Z}})
    }{\int L_g d\bar P_X(\cdot|\mathcal{I}_X)(x^{\mathbb{Z}})}
    \label{PlainInverseAMS}
\end{align}
both for $P_X$-almost all and $\bar P_X$-almost all
$x^{\mathbb{Z}}$. Thus $P_Y$ is AMS.

Suppose next that $P_X$ is given by (\ref{Plain}) for an AMS measure
$P_Y$ with $\bar P_Y(g(\Omega_X)|\mathcal{I}_Y)>0$ $P_Y$-almost
surely. Then $\Phi_X(A)$ exists $P_X$-almost surely if
$\Phi_X(A)\circ g^{-1}$ exists $P_Y$-almost surely. By the Birkhoff
ergodic theorem \cite{Birkhoff32, Garsia65} and by condition
(\ref{SelfAvoiding}), we obtain
\begin{align}
  \Phi_X(A)(g^{-1}(y^{\mathbb{Z}}))
  &=
    \lim_{n\to\infty}\frac{1}{n}\sum_{i=0}^{n-1}
    \boole{T_X^i g^{-1}(y^{\mathbb{Z}})\in A}
    \nonumber\\
  &=
  \lim_{m\to\infty}\frac{\sum_{j=0}^{m-1} \boole{T_Y^j y^{\mathbb{Z}}\in g(A)}
    }{\sum_{j=0}^{m-1} \boole{T_Y^j y^{\mathbb{Z}}\in g(\Omega_X)}}
    \nonumber\\
  &=
    \frac{\bar P_Y(g(A)|\mathcal{I}_Y)(y^{\mathbb{Z}})
    }{\bar P_Y(g(\Omega_X)|\mathcal{I}_Y)(y^{\mathbb{Z}})}
    \label{PlainAMS}
\end{align}
both for $P_Y$-almost all and $\bar P_Y$-almost all
$y^{\mathbb{Z}}$. Thus $P_X$ is AMS.

\subsubsection{Stationary means}

Suppose that probability measures $P_X$ and $P_Y$ are AMS and linked
by (\ref{PlainInverse})--(\ref{Plain}) with $P_Y(g(\Omega_X))=1$. In
this case, conditions (\ref{TrajectoryInverse})--(\ref{Trajectory})
are satisfied with substitutions $P_X\to\bar P_X$ and $P_Y\to\bar P_Y$
by formulas (\ref{StationaryMeanFrequency}), (\ref{PlainInverseAMS}),
and (\ref{PlainAMS}).  Conditions
$\int L_g d\bar P_X(\cdot|\mathcal{I}_X)<\infty$ $P_X$-almost surely
and $\bar P_Y(g(\Omega_X)|\mathcal{I}_Y)>0$ $P_Y$-almost surely are
equivalent by Theorem \ref{theoEnsemble}(i) and formulas
(\ref{TrajectoryInverse})--(\ref{Trajectory}) with substitutions
$P_X\to\bar P_X$ and $P_Y\to\bar P_Y$.

\section{Recurrence times}
\label{secRecurrence}

In this section, we will exhibit an application of the normalized
transport to the problem of successive recurrence times.  Our direct
inspiration for the discovery of mappings (\ref{Ensemble}) and
(\ref{Trajectory}) was learning in the beginning of 2023 that
successive recurrence times are stationary with respect to the
conditional measure. Namely, we have this statement:
\begin{proposition}[cf.\ \mbox{\cite{ChenMoy59}}]
  \label{theoRecurrence}
  Consider a~measurable space $(\Omega,\mathcal{J})$ with an
  automorphism $T:\Omega\to\Omega$ and a~probability measure $P$
  stationary ergodic with respect to $T$. Let $C\in\mathcal{J}$ be
  a~certain event with $P(C)>0$. Define successive recurrence times as
  partial functions
  \begin{align}
    \label{Recurrence}
    R^C_k(\omega)
    :=
    \begin{cases}
      \sup\klam{i\le -1: T^{R^C_{k+1}(\omega)+i}  \omega\in C},
      &k\in-\mathbb{N},
      \\
      0,
      &k=0 \text{ and } \omega\in C,
      \\
      \inf\klam{i\ge 1: T^{R^C_{k-1}(\omega)+i} \omega\in C},
      &k\in\mathbb{N}.
    \end{cases}
  \end{align}
  Process $(R^C_k)_{k\in\mathbb{Z}}$ is stationary ergodic with
  respect to measure $P(\cdot|C)$.
\end{proposition}
\emph{Remark:} In work \cite{ChenMoy59}, stationarity of process
$(R^C_k)_{k\in\mathbb{Z}}$ restricted to $k\in\mathbb{N}$ was
demonstrated by a~direct calculation. The remaining problem whether
process $(R^C_k)_{k\in\mathbb{Z}}$ is ergodic was posted by us in
February 2023 at Mathematics Stack
Exchange.\footnote{\url{https://math.stackexchange.com/questions/4647725/are-successive-recurrence-times-ergodic-for-an-ergodic-process}}
Proposition \ref{theoRecurrence} states that the answer is positive.
This result is an easy corollary of Theorem \ref{theoEnsemble}.
\begin{proof}
We construct process $(Y^C_i)_{i\in\mathbb{Z}}$ with variables
$Y^C_i(\omega):=\boole{T^i\omega\in C}$. Consider measure
$P_Y(B):=P((Y^C_i)_{i\in\mathbb{Z}}\in B)$ on space
$(\Omega_Y,\mathcal{J}_Y)$ with the binary alphabet
$\mathbb{Y}=\klam{\mathbf{0},\mathbf{1}}$.  Next, we consider a~unary
comma-separated code $f:\mathbb{X}\to\mathbb{Y}^*$ with input alphabet
$\mathbb{X}=\mathbb{N}\cup\klam{0}$ and values
\begin{align}
  f(x):=\mathbf{0}^x\mathbf{1}.
\end{align}
Extension $g=f^{\mathbb{Z}}$ is a~self-avoiding code. Using random
variables (\ref{Projections}), we notice equality
\begin{align}
  \label{RecurrenceNormalization}
  P_Y(g(\Omega_X))=P_Y(Y_0=\mathbf{1})=P(C).
\end{align}
We define random variables $R_k$ via (\ref{RecurrenceY}) with
$B=(Y_0=\mathbf{1})$ on $(\Omega_Y,\mathcal{J}_Y)$ and random
variables $L_k:=\abs{f(X_k)}$ on $(\Omega_X,\mathcal{J}_X)$. We
observe that $S_k=R_k(g(X^{\mathbb{Z}}))$ and
$R_k((Y^C_i)_{i\in\mathbb{Z}})=R^C_k$.  Applying the normalized
transport (\ref{Ensemble})--(\ref{EnsembleInverse}) and observation
(\ref{RecurrenceNormalization}), we obtain
\begin{align}
  \label{RecurrenceTransport}
  P_X((L_k)_{k\in\mathbb{Z}}\in B)
  =
  P_Y((R_k)_{k\in\mathbb{Z}}\in B|Y_0=\mathbf{1})
  =
  P((R^C_k)_{k\in\mathbb{Z}}\in B|C)
  .
\end{align}
Recall that measure $P$ is stationary ergodic with respect to $T$.
This implies that measure $P_Y$ is stationary ergodic with respect to
$T_Y$.  Hence measure $P_X$ is stationary ergodic with respect to
$T_X$ by Theorem \ref{theoEnsemble}. This implies that process
$(L_k)_{k\in\mathbb{Z}}$ is stationary ergodic with respect to
$P_X$. Hence process $(R^C_k)_{k\in\mathbb{Z}}$ is stationary ergodic
with respect to measure $P(\cdot|C)$ by (\ref{RecurrenceTransport}).  
\end{proof}

\section{Entropy rate}
\label{secEntropy}


It is natural to ask how some information measures are related for two
processes connected by the normalized transport.  The Shannon block
entropy of a process $(X_i)_{i\in\mathbb{Z}}$ with respect to a
stationary measure $P_X$ is
\begin{align}
  H_X(n):=-\sum_{x_1^n\in\mathbb{X}^n} P_X(X_1^n=x_1^n)\log P_X(X_1^n=x_1^n),
\end{align}
where $x_j^k=(x_j,x_{j+1},...,x_k)$ denotes a block of consecutive
symbols. Similarly, we have
\begin{align}
  H_Y(n):=-\sum_{y_1^n\in\mathbb{Y}^n} P_Y(Y_1^n=y_1^n)\log P_Y(Y_1^n=y_1^n)
\end{align}
for a process $(Y_i)_{i\in\mathbb{Z}}$ and a stationary measure $P_Y$.
The specific question is to characterize how entropies $H_X(n)$ and
$H_Y(n)$ are related if measures $P_X$ and $P_Y$ are connected by the
normalized transport.

This sort of problems was investigated in the past for more specific
mappings between $P_X$ and $P_Y$, cf. \cite[Example 6]{GrayKieffer80},
\cite[Corollary 4.2.5]{Gray90}, \cite{TimoBlackmoreHanlen2007},
\cite[Proposition 6.4]{Debowski10}, \cite[Chapters 10 and
11]{Debowski21}.  Relating non-asymptotic quantities $H_X(n)$ and
$H_Y(n)$ may be cumbersome but things become simpler in the limit. Let
us consider the entropy rates
\begin{align}
  h_X&:=\lim_{n\to\infty}\frac{H_X(n)}{n},
  &
    h_Y&:=\lim_{n\to\infty}\frac{H_Y(n)}{n}.
\end{align}

We state the following result that generalizes \cite[Proposition
6.4]{Debowski10}, which was originally stated for extensions of
comma-separated codes.
\begin{proposition}
  Let stationary ergodic measures $P_X$ and $P_Y$ be connected by the
  normalized transport (\ref{EnsembleInverse})--(\ref{Ensemble}) with
  $\int L_g dP_X<\infty$.  Then we have
  \begin{align}
    \label{HYHX}
    h_Y=\frac{h_X}{\int L_g dP_X}.
  \end{align}
\end{proposition}
\emph{Remark:} Since the entropy rate of a non-ergodic measure is the
expectation of the entropy rate of its random ergodic component
\cite{GrayDavisson74, GrayDavisson74b} then formula (\ref{HYHX}) holds
also for stationary non-ergodic measures $P_X$ and $P_Y$ connected by
the generalized normalized transport
(\ref{TrajectoryInverse})--(\ref{Trajectory}) if the conditional
expectation $\int L_g dP_X(\cdot|\mathcal{I}_X)$ is constant.
\begin{proof}
  Let $s:\Omega_X\to\Omega_W$ be the chunked code induced by the
  self-avoiding code $g:\Omega_X\to\Omega_W$. Consider process
  $(W_i)_{i\in\mathbb{Z}}:=s((X_i)_{i\in\mathbb{Z}})$.  According to
  \cite[Lemma 4.2.3]{Gray90} and \cite[Corollary 4.2.5]{Gray90}, a
  stationary code conserves the ergodicity and stationarity and
  decreases the entropy rate of the coded process.  Since any chunked
  code is a bijective stationary code by Proposition
  \ref{theoChunkedCode} then $(W_i)_{i\in\mathbb{Z}}$ is stationary
  ergodic with respect to measure $P_X$ and we have both $h_W\le h_X$
  and $h_X\le h_W$ for the entropy rate $h_W$ of process
  $(W_i)_{i\in\mathbb{Z}}$.  Thus $h_W=h_X$. On the other hand, by the
  Shannon-McMillan-Breiman theorem
  \cite{Shannon48,Breiman57,Chung61,AlgoetCover88} for process
  $(W_i)_{i\in\mathbb{Z}}$,
  \begin{align}
    \lim_{n\to\infty}\frac{-\log P_X(W_1^n)}{n}
    &=h_W \text{ $P_X$-almost surely}.
  \end{align}
  Since $\abs{W_{i+1}}=L_g\circ T^i$ and $(W_i)_{i\in\mathbb{Z}}$ is
  stationary ergodic then by the Birkhoff ergodic theorem,
  \begin{align}
    \lim_{n\to\infty}\frac{1}{n}\sum_{i=1}^n\abs{W_i}
    &=\int L_g dP_X \text{ $P_X$-almost surely}.
  \end{align}
  Hence we derive
  \begin{align}
    \label{WSMB}
    \lim_{n\to\infty}\frac{-\log P_X(W_1^n)}{\sum_{i=1}^n\abs{W_i}}
    &=\frac{h_X}{\int L_g dP_X} \text{ $P_X$-almost surely}.
  \end{align}

  Observation (\ref{WSMB}) should be contrasted with the
  Shannon-McMillan-Breiman theorem for process
  $(Y_i)_{i\in\mathbb{Z}}$,
  \begin{align}
    \label{YSMB}
    \lim_{n\to\infty}\frac{-\log P_Y(Y_1^n)}{n}
    &=h_Y \text{ $P_Y$-almost surely}.
  \end{align}
  We have
  \begin{align}
     P_X(A)=\frac{P_Y(g(A))}{P_Y(g(\Omega_X))},
  \end{align}
  where $P_Y(g(\Omega_X))=1/\int L_g dP_X>0$. Hence if $P_X(A)=1$ and
  $P_Y(B)=1$  then
  $g(A)\cap B\neq\emptyset$.
  In particular, we may take events
  \begin{align}
    A
    &:=\okra{\lim_{n\to\infty}\frac{-\log P_X(W_1^n)}{\sum_{i=1}^n\abs{W_i}}
    =\frac{h_X}{\int L_g dP_X}},
    \\
    B
    &:=\okra{\lim_{n\to\infty}\frac{-\log P_Y(Y_1^n)}{n}
      =h_Y}.
  \end{align}
  Let us write $a\sqsubset b$ when string $a$ is a prefix of sequence
  $b$.  The images of cylinders are
  $g(w_1^n\sqsubset W_1^\infty)=(y_1^l\sqsubset Y_1^\infty)$, where
  string $y_1^l$ is the concatenation of strings $w_1, w_2, ...,
  w_n$. Hence we obtain
  \begin{align}
    g(A)
    &=
      g\okra{\lim_{n\to\infty}\frac{-\log P_X(W_1^n)}{\sum_{i=1}^n\abs{W_i}}
      =\frac{h_X}{\int L_g dP_X}}
      \nonumber\\
    &=
      g\okra{\lim_{n\to\infty}\frac{-\log P_Y\circ g(W_1^n)
      +\log P_Y\circ g(\Omega_X)}{\sum_{i=1}^n\abs{W_i}}
      =\frac{h_X}{\int L_g dP_X}}
      \nonumber\\
    &=
      g\okra{\lim_{n\to\infty}\frac{-\log P_Y\circ g(W_1^n)
      }{\sum_{i=1}^n\abs{W_i}}
      =\frac{h_X}{\int L_g dP_X}}
      \nonumber\\
    &\subset
      \okra{\liminf_{n\to\infty}\abs{\frac{-\log P_Y(Y_1^n)}{n}
      -\frac{h_X}{\int L_g dP_X}}=0}.
  \end{align}
  Consequently,
  \begin{align}
    \emptyset\neq g(A)\cap B\subset
    \okra{h_Y=\lim_{n\to\infty}\frac{-\log P_Y(Y_1^n)}{n}
    =\frac{h_X}{\int L_g dP_X}}.
  \end{align}
   Hence  the claim follows.
\end{proof}

\section{Open problems}
\label{secProblems}

To conclude the article, let us sketch some open problems.  The
classes of self-avoiding codes and bijective stationary codes are
worth further inspection in the future.  We suspect that there must be
some sort of separators or other synchronizing signals in
a~self-avoiding code. Nonetheless, these separators may be placed
quite arbitrarily as Example \ref{exEmbedded} has shown. Potentially,
synchronizing separators may be distributed sparsely, if they are
combined with fix-free codes \cite{GillmanRivest95, Debowski10}. Thus,
it is not immediate how complicated self-avoiding codes can be in
general---and we may allow for some algorithmic randomness in their
definition.  Moreover, it may be natural to seek for generalizations
of mappings $g=g_1^{-1}\circ g_2$, where $g_1$ and $g_2$ are
self-avoiding codes. These may be a natural extension of self-avoiding
codes that enjoys many properties of theirs.

There may be also interesting connections with formal linguistics.  We
suspect that further examples of complex self-avoiding codes may be
inspired by quirks of the orthography of natural language. The
context-dependent writing of the verb negation in English described as
Example \ref{exNot} is only a starting point.  Speaking more
abstractly, interesting examples of self-avoiding codes may be
supplied by idiosyncrasies of translation between texts in natural
language and their various formal representations, including
translation between different human languages---or plausibly also
between programming languages. It may be worth a systematic linguistic
study to classify which types of self-avoiding codes are plausible in
particular languages and translation systems.  Such a detailed
linguistic analysis may exceed simple word orthography, used here as a
teaser.  Certainly, such a study goes beyond the scope of this article
but its conclusions may enrich information theory.

\bibliographystyle{IEEETran}

\bibliography{0-journals-abbrv,0-publishers-abbrv,ai,ql,mine,tcs,books,nlp}

\begin{thebibliography}{10}
\providecommand{\url}[1]{#1}
\csname url@rmstyle\endcsname
\providecommand{\newblock}{\relax}
\providecommand{\bibinfo}[2]{#2}
\providecommand\BIBentrySTDinterwordspacing{\spaceskip=0pt\relax}
\providecommand\BIBentryALTinterwordstretchfactor{4}
\providecommand\BIBentryALTinterwordspacing{\spaceskip=\fontdimen2\font plus
\BIBentryALTinterwordstretchfactor\fontdimen3\font minus
  \fontdimen4\font\relax}
\providecommand\BIBforeignlanguage[2]{{%
\expandafter\ifx\csname l@#1\endcsname\relax
\typeout{** WARNING: IEEEtran.bst: No hyphenation pattern has been}%
\typeout{** loaded for the language `#1'. Using the pattern for}%
\typeout{** the default language instead.}%
\else
\language=\csname l@#1\endcsname
\fi
#2}}

\bibitem{GrayKieffer80}
R.~M. Gray and J.~C. Kieffer, ``Asymptotically mean stationary measures,''
  \emph{Ann.\ Probab.}, vol.~8, pp. 962--973, 1980.

\bibitem{Debowski21}
{\L}.~D\k{e}bowski, \emph{Information Theory Meets Power Laws: Stochastic
  Processes and Language Models}.\hskip 1em plus 0.5em minus 0.4em\relax Wiley
  \& Sons, 2021.

\bibitem{ChenMoy59}
S.-T. {Chen Moy}, ``Successive recurrence times in a stationary process,''
  \emph{Ann.\ Math.\ Statist.}, vol.~30, no.~4, pp. 1254--1257, 1959.

\bibitem{Wolff80}
J.~G. Wolff, ``Language acquisition and the discovery of phrase structure,''
  \emph{Lang.\ Speech}, vol.~23, pp. 255--269, 1980.

\bibitem{DeMarcken96}
C.~G. de~Marcken, ``Unsupervised language acquisition,'' Ph.D. dissertation,
  Massachussetts Institute of Technology, 1996.

\bibitem{KiefferYang00}
J.~C. Kieffer and E.~Yang, ``Grammar-based codes: {A} new class of universal
  lossless source codes,'' \emph{IEEE Trans.\ Inform.\ Theory}, vol.~46, pp.
  737--754, 2000.

\bibitem{CharikarOthers05}
M.~Charikar, E.~Lehman, A.~Lehman, D.~Liu, R.~Panigrahy, M.~Prabhakaran,
  A.~Sahai, and A.~Shelat, ``The smallest grammar problem,'' \emph{IEEE Trans.\
  Inform.\ Theory}, vol.~51, pp. 2554--2576, 2005.

\bibitem{SenrichHaddowBirch16}
R.~Sennrich, B.~Haddow, and A.~Birch, ``Neural machine translation of rare
  words with subword units,'' in \emph{Proceedings of the 54th Annual Meeting
  of the Association for Computational Linguistics (Volume 1: Long Papers)},
  K.~Erk and N.~A. Smith, Eds.\hskip 1em plus 0.5em minus 0.4em\relax
  Association for Computational Linguistics, 2016, pp. 1715--1725.

\bibitem{VieiraOthers24}
T.~Vieira, B.~LeBrun, M.~Giulianelli, J.~L. Gastaldi, B.~DuSell, J.~Terilla,
  T.~J. O'Donnell, and R.~Cotterell, ``From language models over tokens to
  language models over characters,'' 2024,
  \url{https://arxiv.org/abs/2412.03719}.

\bibitem{CariolaroPierobon77}
G.~Cariolaro and G.~Pierobon, ``Stationary symbol sequences from
  variable-length word sequences,'' \emph{IEEE Trans.\ Inform.\ Theory},
  vol.~23, pp. 243--253, 1977.

\bibitem{TimoBlackmoreHanlen2007}
R.~Timo, K.~Blackmore, and L.~Hanlen, ``On the entropy rate of word-valued
  sources,'' in \emph{Proceedings of the Telecommunication Networks and
  Applications Conference, ATNAC 2007}, 2007, pp. 377--382.

\bibitem{Cox62}
D.~R. Cox, \emph{Renewal Theory}.\hskip 1em plus 0.5em minus 0.4em\relax
  London: Methuen and Company, 1962.

\bibitem{Elga00}
A.~Elga, ``Self-locating belief and the {Sleeping Beauty} problem,''
  \emph{Analysis}, vol.~60, pp. 143--147, 2000.

\bibitem{Debowski10}
{\L}.~D\k{e}bowski, ``Variable-length coding of two-sided asymptotically mean
  stationary measures,'' \emph{J.\ Theor.\ Probab.}, vol.~23, pp. 237--256,
  2010.

\bibitem{Debowski11b}
------, ``On the vocabulary of grammar-based codes and the logical consistency
  of texts,'' \emph{IEEE Trans.\ Inform.\ Theory}, vol.~57, pp. 4589--4599,
  2011.

\bibitem{Debowski21b}
------, ``A refutation of finite-state language models through {Zipf's} law for
  factual knowledge,'' \emph{Entropy}, vol.~23, p. 1148, 2021.

\bibitem{Debowski23}
------, ``A simplistic model of neural scaling laws: Multiperiodic {Santa Fe}
  processes,'' 2023, \url{https://arxiv.org/abs/2302.09049}.

\bibitem{Zipf35}
G.~K. Zipf, \emph{The Psycho-Biology of Language: {An} Introduction to Dynamic
  Philology}.\hskip 1em plus 0.5em minus 0.4em\relax Houghton Mifflin, 1935.

\bibitem{KaplanOther20}
J.~Kaplan, S.~McCandlish, T.~Henighan, T.~B. Brown, B.~Chess, R.~Child,
  S.~Gray, A.~Radford, J.~Wu, and D.~Amodei, ``Scaling laws for neural language
  models,'' 2020, \url{https://arxiv.org/abs/2001.08361}.

\bibitem{Hutter21}
M.~Hutter, ``Learning curve theory,'' 2021,
  \url{https://arxiv.org/abs/2102.04074}.

\bibitem{Hilberg90}
W.~Hilberg, ``{Der bekannte Grenzwert der redundanzfreien Information in Texten
  --- eine Fehlinterpretation der Shannonschen Experimente?}'' \emph{Frequenz},
  vol.~44, pp. 243--248, 1990.

\bibitem{CrutchfieldFeldman03}
J.~P. Crutchfield and D.~P. Feldman, ``Regularities unseen, randomness
  observed: {The} entropy convergence hierarchy,'' \emph{Chaos}, vol.~15, pp.
  25--54, 2003.

\bibitem{Debowski15d}
{\L}.~D\k{e}bowski, ``Hilberg exponents: New measures of long memory in the
  process,'' \emph{IEEE Trans.\ Inform.\ Theory}, vol.~61, pp. 5716--5726,
  2015.

\bibitem{Gray90}
R.~M. Gray, \emph{Entropy and Information Theory}.\hskip 1em plus 0.5em minus
  0.4em\relax Springer, 1990.

\bibitem{Jelinek97}
F.~Jelinek, \emph{Statistical Methods for Speech Recognition}.\hskip 1em plus
  0.5em minus 0.4em\relax The MIT Press, 1997.

\bibitem{Gray09}
R.~M. Gray, \emph{Probability, Random Processes, and Ergodic Properties}.\hskip
  1em plus 0.5em minus 0.4em\relax Springer, 2009.

\bibitem{Rechard56}
O.~W. Rechard, ``Invariant measures for many-one transformations,'' \emph{Duke
  Math.\ J.}, vol.~23, pp. 477--488, 1956.

\bibitem{Birkhoff32}
G.~D. Birkhoff, ``Proof of the ergodic theorem,'' \emph{Proc.\ Nat.\ Acad.\
  Sci.\ Uni.\ Stat.\ Amer.}, vol.~17, pp. 656--660, 1932.

\bibitem{Garsia65}
A.~M. Garsia, ``A simple proof of {E. Hopf's} maximal ergodic theorem,''
  \emph{J.\ Math.\ Mech.}, vol.~14, pp. 381--382, 1965.

\bibitem{GrayDavisson74}
R.~M. Gray and L.~D. Davisson, ``Source coding theorems without the ergodic
  assumption,'' \emph{IEEE Trans.\ Inform.\ Theory}, vol.~20, pp. 502--516,
  1974.

\bibitem{GrayDavisson74b}
------, ``The ergodic decomposition of stationary discrete random processses,''
  \emph{IEEE Trans.\ Inform.\ Theory}, vol.~20, pp. 625--636, 1974.

\bibitem{Shannon48}
C.~Shannon, ``A mathematical theory of communication,'' \emph{Bell Syst.\
  Tech.\ J.}, vol.~30, pp. 379--423,623--656, 1948.

\bibitem{Breiman57}
L.~Breiman, ``The individual ergodic theorem of information theory,''
  \emph{Ann.\ Math.\ Statist.}, vol.~28, pp. 809--811, 1957.

\bibitem{Chung61}
K.~L. Chung, ``A note on the ergodic theorem of information theory,''
  \emph{Ann.\ Math.\ Statist.}, vol.~32, pp. 612--614, 1961.

\bibitem{AlgoetCover88}
P.~H. Algoet and T.~M. Cover, ``A sandwich proof of the
  {Shannon-McMillan-Breiman} theorem,'' \emph{Ann.\ Probab.}, vol.~16, pp.
  899--909, 1988.

\bibitem{GillmanRivest95}
D.~Gillman and R.~L. Rivest, ``Complete variable-length ``fix-free'' codes,''
  \emph{Designs\ Cod.\ Cryptogr.}, vol.~5, pp. 109--114, 1995.

\end{thebibliography}

\end{document}